\documentclass[a4paper]{article}
\usepackage{amsmath,amssymb}
\newtheorem{thm}{Theorem}
\newtheorem{prop}{Proposition}
\newtheorem{cor}{Corollary}
\newtheorem{lem}{Lemma}

\newtheorem{re}{Remark}
\newtheorem{ex}{Example}
\newenvironment{proof}{
                        \noindent{\bf\small Proof.}\small}
                                       {\hfill {$\mathbf \Box$}\medskip}
\newcommand{\got}[1]{\mathfrak{#1}}
\newcommand{\g}{\got{g}}
\newcommand{\J}{\got{J}}
\newcommand{\gl}{\got{gl}}
\renewcommand{\r}{\got{r}}
\renewcommand{\a}{\got{a}}
\renewcommand{\b}{\got{b}}
\newcommand{\s}{\got{s}}
\newcommand{\m}{\got{m}}
\newcommand{\n}{\got{n}}
\renewcommand{\u}{\got{u}}
\renewcommand{\k}{\got{k}}
\newcommand{\BL}{\got{B}}
\newcommand{\p}{\got{p}}
\newcommand{\h}{\got{h}}
\newcommand{\M}{\Delta}
\newcommand{\V}{\mathrm{V}}
\newcommand{\K}{\mathbb{K}}
\newcommand{\N}{\mathbb{N}}

\newcommand{\A}{\mathcal{A}}
\newcommand{\B}{\mathcal{B}}
\newcommand{\gi}{\widehat{\g}}
\newcommand{\Ni}{\hat{\n}}
\newcommand{\Cg}{C^\infty\left(\g\right)}
\newcommand{\Der}{\mathrm{Der}\,}
\newcommand{\ad}{\mathrm{ad}\,}
\newcommand{\G}{\Gamma}
\renewcommand{\d}{\delta}
\makeatletter
\@addtoreset{equation}{section}

\makeatother

\title{Note on The Generalized Derivation Tower Theorem for Lie algebras}
\author{Toukaiddine Petit\footnote{{\tt Author supported by the Scientific Programme NOG of the European Science Foundation}. 
}\\e-mail: toukaiddine.petit@ua.ac.be\\ Freddy Van Oystaeyen\footnote{{\tt Acknowledging the EC project Liegrits MCRTN 505078}. 
}\\e-mail: Fred.vanoystaeyen@ua.ac.be
\\Department Wiskunde en Informatica, Universiteit Antwerpen\\
B-2020 (Belgium)}   
\date{}

\begin{document}
\maketitle

\subsection*{Abstract:}
We provide an algorithm for decomposing a finite-dimensional Lie algebra over a field of characteristic $0$ permitting to generalize the derivation tower theorem for Lie algebras, is proved by E. Schenkman \cite{Sc}.
\tableofcontents
\newpage
\section*{Introduction}
A Lie algebra $\g$ over a field $\K$ of characteristic $0$, is called complete if the center of $\g$ is trivial and all derivations of $\g$ are inner. Let $(\g_n)$ be the series of algebras defined by $\g_0=\g$ and $\g_{n+1}=\Der\g_n$, it is called the tower of derivation algebras of $\g$. E. Schenkman proved that the derivation tower theorem which asserts that if the center of $\g$ is trivial, then the derivation tower of $\g$ terminates with a complete Lie algebras $\gi$, \cite{Sc}. In this note we revisite this theorem aiming to provide an explicite construction of the limit Lie algebra $\gi$. The method is based on $\Gamma$-decomposition in terms of so-called $\Gamma$-triples which are essentially unique, cf Theorem \ref{n5}. 2. This allows to characterize completeness of $\g$ in terms of the representation $\mu$ associated to a $\Gamma$-triple, cf Theorem \ref{n9}. 3. Both theorems cited providece us with a technique of decomposing Lie algebras that allows to construct the limit $\gi$ when it exists: In Section 2 we carry out this explicite construction under the assomption that the center of $\g_n$ is trivial, i.e. the center of $\Der(\g_{n+1})$ is zero. The form of $\gi$ given in (\ref{E2.11}) follows from the main result Theorem \ref{n10} concerning Lie algebras with trivial center.\\
In Section 3 the general case is considered allowing the bad case when the sequence of the dimension increases divergently. The remaining cases are classified in two classes, the first case dealt with by Theorem \ref{n10} and a second class allowing to discribe $\gi$ as $\K\times\left[\gi,\gi\right]$, the latter being a perfect complete Lie algebra. So, excluding the bad case we obtain concrete structure results for the limit Lie algebra $\gi$, cf Theorem \ref{n17}.
\section{$\Gamma$-Decomposition}
Throughout this paper, $\g$ is a finite-dimensional Lie algebra over a field $\K$ of characteristic $0$, $Z\left(\g\right)$ its center, $\r$ its radical, $\n$ its largest nilpotent ideal, $\Cg$ is the intersection of the ideals $C^p\left(\g\right)$ of the central descending sequence of $\g$ and $\Der\g$ its Lie algebra of derivations. The Lie algebra $\Der\g$ is then algebraic \cite[p. 179]{Ch}. The Lie algebra $\ad\g$ is an ideal of  $\Der\g$. Let $e\left(\ad\g\right)$ be the smallest Lie algebra which is algebraic in $\g$ and contains $\ad\g$ \cite[p. 173]{Ch}. Then we have 
\begin{equation}
	\ad\g\subset e\left(\ad\g\right)\subset\Der\g.
\end{equation}
A Lie subalgebra $\Gamma$ of $\gl\left(\g\right)$ is said to be completely reducible or c.r. if its natural action on $\g$ is semi-simple. This means that this Lie algebra is reductive and its center consists of linear maps which are all semi-simple. A Lie subalgebra of $\gl\left(\g\right)$ is said to be  maximal completely reducible or m.c.r. if it is maximal among the c.r. Lie subalgebras. Two m.c.r. Lie subalgebras of $\gl\left(\g\right)$ are isomorphic \cite{Mo}.
Let $u\in\gl\left(\g\right)$ and let $u=u\mid_{\mathrm{S}}+u\mid_{\mathrm{N}}$ be its Jordan decomposition with $\u\mid_{\mathrm{S}}$ ($u\mid_{\mathrm{N}}$ resp. ) its semi-simple (nilpotent resp.) component. If $\u\subset\gl\left(\g\right)$ is a subspace, we will denote by $\u\mid_{\mathrm{S}}$ ($\u\mid_{\mathrm{N}}$ resp.) the set of semi-simple (nilpotent resp.) components of the Jordan decomposition.
\newpage
\begin{lem}\label{n1}Let $\g$ be a Lie algebra over $\K$ and $\G$ be a c.r Lie subalgebra of $\Der\g$. Then $\g$ satisfies
\begin{enumerate}
	\item $\g=\g^\G\oplus\G\cdot\g,  \left[\g^\G,\G\cdot\g\right]\subset\G\cdot\g,$
where $\g^\G:=\left\{x\in\g:\alpha\cdot x=0,\forall \alpha\in\G\right\},\\\G\cdot\g:=\left\{\alpha\cdot\g,\forall\alpha\in\Gamma\right\}$,
	\item If we set $\p:=\G\cdot\g+\left[\G\cdot\g,\G\cdot\g\right]$, then $\p$ is an ideal of $\g$ generated by $\G\cdot\g$ such that $C^p\left(\g\right)=\p+C^p\left(\g^\G\right)$ $\forall p\in\N\cup\left\{\infty\right\}$,
	\item There exists a Levi subalgebra $\s$ of $\g$ such that $\G\cdot\s\subset\s$,
	\item If we set $\a:=\left(\ad_\g\right)^{-1}\left(\ad\g\cap\G\right)$ then $\a$ is a reductive Lie subalgebra of $\g$ satisfying $\G\cdot\a\subset\a$.
	\end{enumerate}
\end{lem}
\begin{proof}
The natural action of $\Gamma$ on $\g$ being semi-simple, we then have $\g=\g^\G\oplus\G\cdot\g$. For all $x\in\g^\G$, $y\in\g$ and $\alpha\in\Gamma\subset\Der\g$ then 
\begin{equation}
	\alpha\cdot \left[x,y\right]=\left[\alpha\cdot x,y\right]+\left[x,\alpha\cdot y\right]=\left[x,\alpha\cdot y\right]
\end{equation}
and
\begin{equation}
	\left[\g^\G,\G\cdot\g\right]\subset\G\cdot\g.
\end{equation}
The statement 1. holds. Since $\left[\g^\G,\G\cdot\g\right]\subset\G\cdot\g$ we deduce that $\p$ is an ideal of $\g$ and generated by $\G\cdot\g$. The rest of the statement 2. is obvious. $\Gamma$ being reductive, there exists a Levi subalgebra $\s$ of $\g$ such that
\begin{equation}
	\Gamma:=\ad\s\oplus Z\left(\Gamma\right)\quad\mathrm{and}\quad\left[\Gamma,\Gamma\right]=\ad\s.
\end{equation}
Since $\Gamma$ is not maximal then there exists a m.c.r subalgebra $\Gamma_{max}$ of $\Der\g$ containing $\Gamma$ such that $\Gamma_{max}\cdot\s\subset\s$ and a fortiori $\Gamma\cdot\s\subset\s$.
\end{proof}
\begin{lem}\label{n2}Let $\g$ be a Lie algebra over $\K$. Then there exists a nilpotent Lie subalgebra $\h$ of $\g$ such that $\g=\h+C^\infty\left(\g^\G\right)$.
\end{lem}
\begin{proof}
This is obvious if $\dim\g=0$. We reason by induction on $n:=\dim\g$. We may assume that $\g$ is not nilpotent else we set $\g:=\h$. There exists $x\in\g$ such that $\delta:=\ad x\mid_{S}$ is different from zero and $\g=\delta\cdot\g\oplus\g^{\delta}$. We have $\dim\g^{\delta}<n$ and $\g^\delta=\h+C^{\infty}\left(\g^{\delta}\right)$ by the induction hypothesis. The inclusion $\delta\cdot\g\subset C^\infty\left(\g\right)$ yields $\delta\left(\g\right)\subset C^{\infty}\left(\g^\G\right)$ and $\g=\h+C^\infty\left(\g^\G\right)$.
\end{proof}
\begin{cor}\label{n3}
Let $\g$ be a Lie algebra over $\K$. There exists a Levi subalgebra $\s$ of $\g$ and a nilpotent Lie subalgebra $\h$ of $\r$ such that 
\begin{enumerate}
	\item $\g=\s+\h+\n$ and $\left[\s,\h\right]=0$,
	\item $e\left(\ad\g\right)=\G\oplus \M$ with $\G:=\ad\s+\ad \h\mid_{\mathrm{S}}$, $\M:=\ad \h\mid_{\mathrm{N}}+\ad\n$ and $\G\cdot \h=0$.
	We will say that $\G$ is a m.c.r Lie subalgebra of $e\left(\ad\g\right)$ associated to $\left(\s,\h\right)$.
\end{enumerate}

\end{cor}
\begin{proof}
By Lemma \ref{n2} there exists a nilpotent Lie subalgebra $\h$ of $\r$ such that 
\begin{equation}
	\r^{\ad\s}=\h+C^{\infty}\left(\r^{\ad\s}\right)\subset \h+\n \quad\mathrm{and}\quad\r=\r^{\ad\s}\oplus\left[\s,\r\right].
\end{equation}
We have $\left[\s,\r\right]\subset\n$, hence $\g=\s+\h+\n$ and $\left[\s,\h\right]=0$ and $\h$ is nilpotent.
It follows that 
\begin{equation}
	\ad\g=\ad\s+\ad \h+\ad\n\quad\mathrm{and}\quad e\left(\ad\g\right)=\ad\s+e\left(\ad\h \right)+\ad\n.
\end{equation}
The Lie algebra $e\left(\ad \h\right)$ is nilpotent and admits a Chevalley decomposition 
\begin{equation}
	\ad \h\mid_{\mathrm{N}}\oplus\ad \h\mid_{\mathrm{S}}, \quad\cite{Ch}.
\end{equation}
It follows that 
\begin{equation}
	e\left(\ad\g\right)=\G\oplus \M\quad\mathrm{and}\quad\G\cdot \h=0
\end{equation}
 with 
\begin{equation}
	\G:=\ad\s+\ad \h\mid_{\mathrm{S}}\quad\mathrm{and}\quad\M:=\ad \h\mid_{\mathrm{N}}+\ad\n.
\end{equation}
Hence $\G$ is maximal by construction.
\end{proof}

We introduce the notion of a $\G$-decomposition in $iv)$ hereafter.

\begin{thm}\label{n5}
Let $\g$ be a Lie algebra over $\K$ and $\r$ its radical. 
\begin{enumerate}
	\item There is a bijection of the set of c.r.m subalgebras of $e\left(\ad\g\right)$ into the set of sequences of vector spaces $\left(\s,\k,\m\right)$ of $\g$ such that :
		\item [ i)] $\s$ is a Levi subalgebra of $\g$,
		\item [ ii)] $\k$ is an ideal nilpotent of $\g$ such that $\left[\s,\k\right]=0$,
		\item [ iii)] $\m$ is a subspace of $\r$ such that $\r=\m\oplus\k$ and $\left[\s\oplus\k,\m\right]=\m$, and
		\item [ iv)] $\g =\s\oplus\k\oplus\m$.
		Will call $\left(\s,\k,\m\right)$ a $\G$-triple and $\s\oplus\k\oplus\m$ a $\G$-decomposition of $\g$.
		
	\item Let $\left(\s_i,\k_i,\m_i\right)$ be two $\G_{i}$-triple of $\g$ with $i=1,2$, then there exists an inner automorphism $\gamma $ of $\g$ such that 
\begin{equation}
\G_{2}=\gamma\circ\G_{1}\circ\gamma^{-1},\quad\s_2=\gamma\left(\s_1\right),\quad\m_2=\gamma\left(\m_1\right),\quad\k_2=\gamma\left(\k_1\right).
\end{equation}

\end{enumerate}
\end{thm}
\begin{proof}
If we assume that there exists a such decomposition $\g =\s\oplus\k\oplus\m$, we construct $\G$ by setting 
\begin{equation}
	\G:=\ad\s\oplus\ad \k\mid_{\mathrm{S}}
\end{equation}
i.e. Corollary \ref{n3}.2. Conversely, let $\G$ be a m.c.r Lie subalgebra of $e\left(\ad\g\right)$, we set 
\begin{equation}
	\k:=\g^{\G}\quad\m:=\G\cdot\r\quad\mathrm{and}\quad\s:=\left(\ad_\g\right)^{-1}\left(\ad\g\cap\G\right).
\end{equation}
It is obvious that $\s$ is a Levi subalgebra of $\g$. We have 
\begin{equation}
	\G\cdot\g=\G\cdot\m\oplus\G\cdot\k\oplus\G\cdot\s=\m\oplus\s
\end{equation} 
since 
\begin{equation}
	\G\cdot\m=\m,\quad\G\cdot\s=\s\quad\mathrm{and}\quad\G\cdot\k=0.
\end{equation}
By Lemma \ref{n1}, we have 
\begin{equation}
	\g=\g^\G\oplus\G\cdot\g=\s\oplus\k\oplus\m.
\end{equation}
Since 
\begin{equation}
	\left[\G,\ad\k\right]=\ad\left(\G\cdot\k\right)=0
\end{equation}
it follows that 
\begin{equation}
	\left[\G,\ad\k\mid_{\mathrm{S}}\right]=\left[\G,\ad\k\mid_{\mathrm{N}}\right]=0.
\end{equation}
Hence 
\begin{equation}
	\ad\k\mid_{\mathrm{S}}\subset Z\left(\G\right)\subset\G
\end{equation} 
because $\G$ is m.c.r Lie subalgebra. The Lie algebra 
\begin{equation}
	Z\left(\G\right)+\ad\k=Z\left(\G\right)\oplus\ad\k\mid_{\mathrm{N}}
\end{equation}
is thus nilpotent. Then $\ad\k$ is nilpotent and we conclude that $\k$ is nilpotent.
We have $[\s\oplus\k,\m]\subset\m$ since $\m\subset\r.$
If $[\s\oplus\k,\m]\neq\m$ 
then 
\begin{equation}
	(\ad\s+e(\ad\k))\m\neq\m\quad\mathrm{and}\quad \G\cdot\m\neq\m.
\end{equation}
This is a contradiction. The statement 2. is a consequence of G. D. Mostow's theorem \cite{Mo} applied to $e\left(\ad\g\right)$.
\end{proof}
\begin{cor}\label{n6}
Let $\left(\s,\k,\m\right)$ be a $\G$-triple of $\g$ with $\G$ a c.r.m Lie subalgebra $e\left(\ad\g\right)$. Then $$\Der\g=\ad\s\oplus\ad\m\oplus\left(\Der\g\right)^\G$$ 
where $\left(\Der\g\right)^\G$ is the centralizer of $\G$ in $\Der\g$.
\end{cor}
\begin{proof}
The adjoint representation of $\G$ in $\Der\g$ being semi-simple, then
\begin{equation}
	\Der\g=\left[\Der\g,\G\right]\oplus\left(\Der\g\right)^\G.
\end{equation}
Since 
\begin{equation}
	\left[\G,\Der\g\right]\subset\left[e\left(\ad\g\right),\Der\g\right]\subset\ad\g
\end{equation} 
and 
\begin{equation}
	\left(\Der\g\right)^\G\cap\ad\g=\ad\k
\end{equation}
hence 
\begin{equation}
	\Der\g=\ad\s\oplus\ad\m\oplus\left(\Der\g\right)^\G.
\end{equation}
\end{proof}

\begin{re} 
The vector space $\m$ is not a subalgebra of $\g$ in general. If we consider the solvable Lie algebra generated by $\left\{x_1,x_2,x_3,x_4,x_5\right\}$ with its multiplication defined by
\begin{equation}
	\left[x_1,x_2\right]=x_5,\quad\left[x_1,x_3\right]=x_3,\quad\left[x_1,x_4\right]=-x_4,\quad\left[x_3,x_4\right]=x_5.
\end{equation}
Then 
\begin{equation}
	\G=\K\cdot\ad x_1\mid_{\mathrm{S}},\quad\k=\K\cdot x_1+\K\cdot x_2+\K\cdot x_5,\quad\m=\K\cdot x_3+\K\cdot x_4.
\end{equation}
We observe that $\m$ is not an ideal.
\end{re}

\begin{prop}Let $\left(\s_i,\k_i,\m_i\right)$ be a $\G_i$-triple of the Lie algebra $\g_i$ with each $\G_i$ a c.r.m Lie subalgebras of $e\left(\ad\g_i\right)$ with $i=1,2$. Then $\left(\s_1\times\s_2,\k_1\times\k_2,\m_1\times\m_2\right)$ is a $\G_1\times\G_2$-triple of the Lie algebra $\g_1\times\g_2$.
\end{prop}
\begin{prop}Let $f:\g\rightarrow\g'$ be an Lie algebra epimorphism and $\G$ a c.r.m subalgebra of $e\left(\ad\g\right)$. Then 
\begin{enumerate}
	\item The subalgebra $\G'$ induced by $\G$ on the quotient $\g'$ satisfying $\G'\circ f=f\circ\G$ is a c.r.m subalgebra of $e\left(\ad\g'\right)$.
	\item If $\left(\s,\k,\m\right)$ is a $\G$-triple of $\g$ with $\G$ a c.r.m Lie subalgebra of $e\left(\ad\g\right)$, then $\left(f\left(\s\right),f\left(\k\right),f\left(\m\right)\right)$ is a $\G'$-triple of $\g'$.
\end{enumerate}
\end{prop}
\begin{proof}
It is easy to check that the property 1. We have 
\begin{equation}
	\g'=f\left(\s\right)\oplus f\left(\k\right)\oplus f\left(\m\right)
\end{equation}
where $f\left(\k\right)\oplus f\left(\m\right)$ is the radical of $\g'$. Since 
\begin{equation}
	\left[f\left(\s\right)\oplus f\left(\m\right),f\left(\m\right)\right]=f\left(\m\right),
\end{equation}
then $\left(f\left(\s\right),f\left(\k\right),f\left(\m\right)\right)$ is a $\G'$-triple of $\g'$.
\end{proof}

\begin{lem}\label{n7}
Let $\left(\s,\k,\m\right)$ be a $\G$-triple of $\g$ with $\G$ a c.r.m Lie subalgebra of $e\left(\ad\g\right)$. The inclusion $\left[\k,\m\right]\subset\m$ defines a representation 
$$\mu:\k\rightarrow\gl\left(\m\right),\quad x\mapsto\ad x\mid_{m}.$$
We then have :
\begin{enumerate}
	\item $Z\left(\g\right)=Z\left(\k\right)\cap\ker\mu$.
	\item $\mu$ is injective if and only if $Z\left(\g\right)=0$.
	\item The Lie subalgebra $\Ni$ of $\g$ generated by $\m$ is a nilpotent ideal of $\g$ such that 
$$\Ni=C^\infty\left(\g\right)\cap\n=\m+\left[\m,\m\right],\quad C^p\left(\g\right)=\s+C^p\left(\k\right)+\Ni, \quad\forall p\in\N.$$
\end{enumerate}
\end{lem}

\begin{proof} Let $x\in Z\left(\k\right)\cap\ker\mu$ then 

\begin{equation}
	\left[\\s\oplus\k\oplus\m,x\right]=\left[\k,x\right]=0
\end{equation}
and $x\in Z\left(g\right)$. Conversely, if $x\in Z\left(\g\right)$ we have $\left[\g,x\right]=\left\{0\right\}$, hence 
\begin{equation}
	\G\cdot x\subset\k\quad\mathrm{and}\quad\G^2\cdot x=0.
\end{equation} 
This means that $x\in\k\quad\mathrm{so}\quad x\in Z\left(\k\right)\cap\ker\mu.$
 Hence statement 1 holds. The assumption that $\mu$ is not injective is equivalent to 
\begin{equation}
	Z\left(\g\right)=Z\left(\k\right)\cap\ker\mu\neq 0
\end{equation}
since every non-null nilpotent ideal intersects the center. Hence statement $2$ holds. Using Lemma \ref{n1}.2 and the $\G$-decomposition of $\g$, we deduce statement 3.
\end{proof}

\begin{lem}\label{n8}
Let $\V$ be a vector space over $\K$. Let $\b$,$\a$ be two Lie subalgebras of $\gl\left(\V\right)$ such that $\a\subset\b$. Let $\Phi$ be a map of $\V^2$ into $\g:=\a\oplus\V$. We define a bracket $\left[,\right]$ on $\g$ by
 \begin{eqnarray*} 
 & & \left[x,v\right]=x\cdot v \\
 & &\left[v_1,v_2\right]=\Phi(v_1,v_2)\\
 &&\left[x,y\right]=x\cdot y-y\cdot x
\end{eqnarray*}
with $x,y\in\a$, $v_1,v_2\in\V$.
If this bracket defines a Lie structure on $\g$, then it can be lifted to $\g':=\b\oplus\V$ by setting : 
$$\left[x,\Phi(v_1,v_2)\right]=\Phi\left(x\cdot v_1,v_2\right)+\Phi\left(v,x\cdot v_2\right)$$
for all $\left(x,v_1,v_2\right)\in\b\times\V^2$.
\end{lem}

Let $\left(\s,\k,\m\right)$ be a $\G$-triple of $\g$ with $\G$ a m.c.r Lie subalgebra $e\left(\ad\g\right)$. Let $\delta\in\left(\Der\g\right)^\G$. Then it is easy to check that 
\begin{equation}
	\delta\left(\m\right)=\m,\quad\delta\left(\k\right)=\k
\end{equation} 
and 
\begin{equation}
	\delta\left(\s\right)=0.
\end{equation}
This means that $\delta$ defines a linear map 
\begin{equation}
	\Theta:\left(\Der\g\right)^\G\rightarrow\gl\left(\m\right)\quad\delta\mapsto\delta\mid_{\m}. 
\end{equation} 
The set of derivations $\left(\Der\Ni\right)^\G$ of $\Ni$ which commute with the restriction $\G\mid_{\Ni}$ of $\G$ to $\Ni$, stabilizes $\k\cap\Ni$ and $\m$. There exists an isomorphism of Lie algebras of $\left(\Der\Ni\right)^\G$ into $\left(\Der\Ni\right)^\G\mid_{\m}$. The image of $\Theta$ is contained in $\BL:=\left(\Der\Ni\right)^\G\mid_{\m}$.

Lemma \ref{n7} and Lemma \ref{n8}, entail the following:

\begin{thm}\label{n9}
Let $\left(\s,\k,\m\right)$ be a $\G$-triple of $\g$ with $\G$ a c.r.m subalgebra of $e\left(\ad\g\right)$. Let $\mu:\k\rightarrow\gl\left(\m\right)$ be the representation defined by $\mu\left(x\right)=\ad x\mid_{m}$. Assume $\mu$ is injective, then:
\begin{enumerate}
	\item The map$\quad\Theta$ defines an isomorphism of $\Der\left(\g\right)^\G$ into the normalizer $\mathrm{N}_\B\left(\mu\left(\k\right)\right)$ of $\mu\left(\k\right)$ in $\BL$.
	\item We may identify the Lie algebra $\Der\g$ with the Lie algebra $s\oplus \mathrm{N}_\B\left(\mu\left(\k\right)\right)\oplus\m$ with its law $\Phi$ defined from the bracket $\left[,\right]$ of $\g$ in the following way, for all $\left(x_1,y_1,z_1\right),\left(x_2,y_2,z_2\right)\in\s\times \mathrm{N}_\BL\left(\mu\left(\k\right)\right)\times\m $:
 \begin{eqnarray*} & & \Phi\left(x_1,x_2+y_2+z_2\right)=\left[x_1,x_2\right]+\left[x_1,z_2\right]\\
& & \Phi\left(y_1,y_2+z_2\right)=y_1\cdot y_2-y_2\cdot y_1+y_1\cdot z_2\\
& & \Phi\left(z_1,z_2\right)=\mu(k)+m
\end{eqnarray*}
with $k$ and $m$ the projections of $\left[z_1,z_2\right]$ into $\k$ and $\m$, respectively.
	\item In particular $\g$ is complete if and only if $\mu \left(\k\right)$ is equal to its normalizer in $\BL$.
	\end{enumerate}
\end{thm}
\begin{proof} Let $\d\in\left(\Der\g\right)^\G$ such that $\Theta(\d)=0$, hence:
\begin{equation}
	0=\d\left([x,y]\right)=[\d \left(x\right),y]+[x,\d \left(y\right)]=\mu(\d \left(x\right))\cdot y\quad \forall\quad (x,y)\in\k\times\m.
\end{equation}
It follows that 
\begin{equation}
	\d \left(x\right)\in\ker\mu=0,\quad\mathrm{and}\quad\d\left(\g\right)=\d\left(\k\right)=0.
\end{equation}
 Then injectivity of $\Theta$ follows. Let $\d_1$ be an element of $\mathrm{N}_{\B}(\mu(\k))$. For all $x_1\in\k$, there exits $x_2=\mu^{-1}\left([\d_1,\mu\left(x_1\right)]\right)$ 
since $\mu$ is injective. It follows that 
\begin{equation}
	\left[\d_1,\mu(x_1)\right]=\mu\left(x_2\right),
\end{equation}
so we may define a derivation $\d_2$ of $\k$ by $\d_2\left(x_1\right)=x_2$. \\
The surjectivity derives from the linear application $\d$ defined by 
\begin{equation}
	\d\mid_{\s}=0,\quad\d\mid_{\k}=\d_2,\quad\d\mid_{\m}=\d_1
\end{equation} 
which belongs to $(\Der\g)^{\G}$. Let $\d_3$ be the derivation of $\Ni$ which extends $\d$. Then $\d_3$ is equal to $\d_2$ on $\k\cap\Ni$. Indeed for $x\in\k\cap\Ni$ we have 
\begin{equation}
	\left[\d_3,\ad_{\Ni}\left(x\right)\right]=\ad_{\Ni}(\d_3\left(x\right))
\end{equation}
and by restriction to $\m$ this yieldss 
\begin{equation}
	[\d_1,\mu(x)]=\mu(\d_3\left(x\right))
\end{equation}
which is also equal to $\mu(\d_2\left(x\right))$ by definition of $\d_2$. Then $\d_3(x)=\d_2(x)$
because of the injectivity of $\mu$. \\
It follows that $\d_1\in(\Der\g)^{\G}$ since the equality 
\begin{equation}
	\d_1\left(\left[x,y\right]\right)=[\d_2 (x),y]+\left[x,\d_1\left(y\right)\right] \quad \forall\quad  (x,y)\in\k\times\m
\end{equation}
is equivalent to 
\begin{equation}
	\left[\d_1,\mu\left(x\right)\right]=\mu(\d_2(x)).
\end{equation}
Thus statement 1 holds. By Corollary \ref{n6}, Lemma \ref{n8} and statement 1, we arrive at statement 2.
The statement 3 follows from statements 1 and 2.
\end{proof}

Theorem \ref{n9} allows to identify the Lie algebra $\g$ with the Lie subalgebra $\s\oplus\mu(\k)\oplus\m$ of $\s\oplus\mathrm{N}_{\BL}(\mu(\k))\oplus\m$ via the isomorphism  defined by
\begin{equation}
	s+k+m\mapsto s+\mu(k)+m\quad\forall\quad (s,k,m)\in\s\times\k\times\m.
\end{equation}
We define a Lie algebra $\s\oplus\BL\oplus\m$ with the same law $\Phi$ by taking $y_1$ and $y_2$ in $\BL$ instead of $\mathrm{N}_{\BL}(\mu(\k))$. We verify that this is really a Lie algebra by using Lemma \ref{n8}. If $\g$ is identified with $\s\oplus\mu(\k)\oplus\m$, we have the Lie algebra inclusions:
\begin{equation} 
\g \subset \s\oplus\mathrm{N}_{\BL}(\mu(\k))\oplus\m \subset \s\oplus\BL\oplus\m
\end{equation}
\begin{cor}Let $\left(\s,\k,\m\right)$ be a $\G$-triple of $\g$ with $\G$ a c.r.m subalgebra of $e\left(\ad\g\right)$. The Lie algebra $\g'=\s\oplus\BL\oplus\m$ is complete.
\end{cor}

\begin{proof}
The triple $\left(\s,\B,\m\right)$ of $\g'$ is a $\G'$-triple with $\G'$ obtained by the extension of $\G$  to $\g'$ which is trivial on $\BL$. Using Theorem \ref{n9}.3 for the natural representation $\mu$ of $\BL$ in $\m$, we deduce that $\s\oplus\B\oplus\m$ is complete.
\end{proof}

\section{Derivation tower of Lie algebras: case with trivial center}
We recall that the derivation tower of a Lie algebra $\g$ is the sequence of Lie algebras $\left(\g_n\right)_{n\in\N}$ such that 
\begin{equation}
	\g_0=\g,\quad \mathrm{and}\quad\g_{n+1}=\Der\left(\g_n\right).
\end{equation}
If the center of $\g_n$ is trivial, hence the center of $\Der\left(\g_{n+1}\right)$ is zero. Then we can identify $\g_n$ with $\ad\g_n$, and we have a sequence of ideals:
\begin{equation} 
\g_n\lhd\g_{n+1}\lhd\ldots\lhd\g_{n+k}\lhd\ldots .
\end{equation}
In this case, Schenkman proved that this sequence has a limit $\gi$, \cite{Sc}.
If $\A$ is a Lie subalgebra of a Lie algebra $\BL$, we consider the sequence of normalizers in $\BL$:
\begin{equation}
	\mathrm{N}_\BL^{p+1}\A=\mathrm{N}_\BL(\mathrm{N}_\BL^p\A),\quad\mathrm{N}_\BL^0\A=\A.
\end{equation}
The following sequence of ideals
\begin{equation} 
\A\lhd\ldots\lhd \mathrm{N}_\BL^p\A\lhd \mathrm{N}_\BL^{p+1}\A\lhd\ldots\lhd \mathrm{N}_\BL^\infty\A \subset \BL 
\end{equation}
terminates for an integer $p$ since the dimension is finite, say $\mathrm{N}_\BL^{q+1}\A=\mathrm{N}_\BL^q\A$ 
and the Lie algebra $\mathrm{N}_\BL^q\A$ denoted $\mathrm{N}_\BL^\infty\A$ will be equal to its normalizer in $\BL$. With this notation we state: 
\begin{thm}\label{n10} 
Let $\s\oplus\k\oplus\m$ be a $\G$-decomposition of a Lie algebra $\g$ with trivial center and $\G$ being a m.c.r. Lie subalgebra of $e(\ad\g)$. The sequence of normalizers
\[ \mu(\k)\lhd \mathrm{N}_\BL^1(\mu(\k))\lhd\ldots\lhd \mathrm{N}_\B^n(\mu(\k))\lhd\ldots\lhd \widehat{\mathrm{N}_\B}(\mu(\k)) \]
contained in $\BL=(\Der\n)^\G\mid\m$ terminates at $\widehat{\mathrm{N}_\BL}(\mu(\k))=\mathrm{N}_\BL^q(\mu(\k))$ for some integer $q$ such that $\mathrm{N}_\BL^{q+1}(\mu(\k))=\mathrm{N}_\BL^q(\mu(\k))$. The derivation tower of a Lie algebra $\g$ is given by the Lie subalgebras $\g_n=\s\oplus \mathrm{N}_\BL^n(\mu(\k))\oplus\m$ of $\s\oplus\k\oplus\m$ and terminates at $\gi:=\s\oplus \widehat{\mathrm{N}_\BL}(\mu(\k))\oplus\m$, thus $\gi$ is complete.
\end{thm}

\begin{proof} 
For $\dim\g=0$, $\g$ always admits a decomposition $\s\oplus\k\oplus\m$ associated with $\G$, cf. Theorem \ref {n5}. We reason by induction on $n:=\dim\g$. Let us assume that for an integer $n\geq 0$ we have 
\begin{equation}
	\g_n=\s\oplus \mathrm{N}_\BL^n\left(\mu\left(\k\right)\right)\oplus\m
\end{equation}
associated with the m.c.r. algebra $\G_n$ obtained by the extension of $\G$ to $\g_n$ which is trivial on $\mathrm{N}_\BL^n(\mu(\k))$ using the inclusion given by (1). The Lie algebra $\G_n$ satisfies 
\begin{equation}
	\G_{n}\cdot\r=\m
\end{equation} 
and $\mu_n$ is injective since it is given by the natural representation of $\BL$ in $\m$. Because of Theorem \ref{n9}, it is possible to identify $\g_{n+1}$ with 
\begin{equation}
	\s\oplus \mathrm{N}_\B^{n+1}(\mu(\k))\oplus\m.
\end{equation}
The Lie algebra $\G_n$ acts by the adjoint representation on $\g_{n+1}$ which is equivalent to the extension of $\G_{n}$ to $\g_{n+1}$ trivial on $\mathrm{N}_\BL^{n+1}(\mu(\k))$. We then set 
\begin{equation}
	\G_{n+1}:=\ad\G_n
\end{equation}
and we have 
\begin{equation}
	\G_{n+1}\mid\g=\G_n\mid\g=\G.
\end{equation}
The rest of the proof is obvious.
\end{proof}\\
We then have 
\begin{equation}
	\widehat{\mathrm{N}_\BL}(\mu(\k))=(\widehat{\mathrm{N}}_{(\Der\n)^\G}(\ad\k))\mid_\m
\end{equation}
and the Lie algebra $\gi$ is also given by
\begin{equation}\label{E2.11} \gi=\s\oplus(\widehat{\mathrm{N}}_{(\Der\n)^\G}(\ad\k))\mid_\m\oplus \m \end{equation}
\begin{ex}\ Let $\a$ be a nilpotent Lie subalgebra of $\gl(\V)$ such that $\a\cdot\V=\V.$ 
Let $\g=\a\oplus V$ be a semi-direct product for the natural representation of $\a$ by the abelian Lie algebra $\V$. This decomposition is associated to $\G$ with 
\begin{equation}
	\k=\a,\quad\m=\V,\quad\G=\ad A
\end{equation}
where $A$ is the set of semi-simple components of elements of $\a$. With notation as in Theorem \ref{n10}, $\BL$ is the centralizer of $A$ in $\gl(\V)$ and we have 
\begin{equation}
	\gi= \widehat{\mathrm{N}_\BL}(\a)\oplus\V.
\end{equation}
\end{ex}
\begin{cor}\label{n11}We redecover the majoration given in \cite{Sc}:
\begin{equation}
	\dim\gi \le \dim(\s\oplus\BL\oplus\m) \le \dim\Der\left(\Cg\right)+\dim Z\left(\Cg\right).
\end{equation}
\end{cor}
\begin{proof} 
Let $\g=\s\oplus\k\oplus\m$ be the decomposition of $\g$ associated with $\G$ with 
\begin{equation}
	\Cg=\s\oplus\n
\end{equation}
cf. Corollary \ref{n6}. $2$. The Lie subalgebra
\begin{equation}
	\ad\left(\Cg\right)+\left(\Der\left(\Cg\right)\right)^\G
\end{equation}
of $\Der\left(\Cg\right)$ where $\left(\Der\left(\Cg\right)\right)^\G$ is the centralizer of the restriction of $\G$ to $\Cg$ in $\Der\left(\Cg\right)$, can be written as 
\begin{equation}
	\ad\s\oplus\left(\Der\left(\Cg\right)\right)^\G\oplus\ad\m.
\end{equation}  
The isomorphisms 
\begin{equation}
	\left(\Der\left(\Cg\right)\right)^\G \cong (\Der\Ni)^\G \cong (\Der\Ni)^\G \mid_\m 
\end{equation}
show that its dimension is always greater than 
\begin{equation}
	\dim\s + \dim(\Der\Ni)^\G + \dim\m + \dim Z(\Cg)
\end{equation}
 which is equal to 
\begin{equation}
	\dim(\s\oplus\BL\oplus\m)-\dim Z(\Cg)
\end{equation}
and the result follows.
\end{proof}

The case where the sequence terminates at the first degree is descibed in the following.
\begin{prop}\label{n12}If $\g$ is a Lie algebra with trivial center, then $\Der\g$ is complete if and only if the ideal $\ad\g$ of $\Der\g$ is characteristic.
\end{prop}
\begin{proof} 
The necessity is trivial. Now if $\ad\g$ is a characteristic ideal of $\Der\g$, it is stable under $\Der\left(\Der\g\right)$ denoted by $\Der^2\left(\g\right)$. The image of the restriction  morphism $\rho$ of $\Der^2\left(\g\right)$ to $\ad\g$ is equal to 
\begin{equation}
	\Der\left(\ad\g\right)\cong \ad\left(\Der\g\right).
\end{equation} 
The kernel $\J$ of $\rho$ has zero intersection with $\ad\left(\Der\g\right)$ and $\Der^2\left(\g\right)$ is the direct sum of the ideal $\J$ and $\ad\left(\Der\g\right)$. We have 
\begin{equation}
	\left[\ad\left(\Der\g\right),\J\right]=-\ad\left(\J\cdot\Der\g\right)]=0
\end{equation}
and $\J\cdot\Der\g$ is zero since the center of $\Der\g$ is zero, so $\J=0$. Thus 
\begin{equation}
	\Der^2\left(\g\right)=\ad\left(\Der\g\right).
\end{equation}
\end{proof}
\begin{cor}\label{n13}If $Z\left(\g\right)=0$ and $[\g,\g]=\g$ then $\Der\g$ is complete.
\end{cor}
\section{The Derivation tower of Lie algebras: general case}
Let $\left(\g_n\right)_{n\in\N}$ be the derivation tower of a Lie algebra $\g$. We now consider the general case.
The ideal $\mathrm{I}$ of $\Der\g$ of derivations which commute with $\ad\g$ is the set of derivations of images contained in $Z\left(\g\right)$. Hence $\mathrm{I}$ vanishes on $\left[\g,\g\right]$ and contains the center of $\Der\g$. If $Z\left(\g\right)$ or $\left[\g,\g\right]=\g$ then $\mathrm{I}=0$. If $\g$ is the direct product $\g_1\times\g_2$, we denote $\mathrm{I}_{ij}$ where $i$ is different to $j$, for the set of derivations vanishing on $\left[\g_i,\g_i\right]\times\g_j$ and their images contained in $Z\left(\g\right)$.
\begin{lem}\label{n15} We set 
\begin{equation}
	\overline{\Der}\g_k=\left\{f\in\Der\left(\g_1\times\g_2\right):f\mid{0\times\g_k=0}\right\}
\end{equation}
for $ k=1,2$. Then 
\begin{equation}
	\Der\left(\g_1\times\g_2\right)=\overline{\Der}\g_1\oplus\overline{\Der}\g_2\oplus \mathrm{I}_{12}\oplus \mathrm{I}_{21}.
\end{equation}
\end{lem}
\begin{proof} We decompose each derivation into four linear maps $\g_i \to \g_j$ for $i,j \in \{1,2\}$ by expressing the derivation property.
\end{proof}
\begin{lem}\label{n16}The sequence $\left(\g_n\right)$ defined by a Lie algebra $\g=\K\times\a$ such that $\dim\left(\Der\g\right)=\dim\g$ belongs to case $1$ or $2$ of Theorem \ref{n14} for $\a=\left[\gi,\gi\right]$.
\end{lem}
\begin{proof} 
Lemma \ref{n15} shows that 
\begin{equation}
	\Der\g=\K\epsilon\oplus\overline{\Der}\a\oplus\mathrm{I}_{12}\oplus \mathrm{I}_{21}
\end{equation}
where $\epsilon$ is the identity on $\K$ and $0$ on $\a$. So 
\begin{equation}
	\dim\left(\Der\g\right) = \dim\g = 1 + \dim\left(\Der\a\right) + \dim Z\left(\a\right) + \dim\left(\a/\left[\a,\a\right]\right). 
\end{equation}
We necessarily have 
\begin{equation}
	\left[\a,\a\right]=\a\quad\mathrm{and}\quad\Der\a=\ad\a.
\end{equation}
If $Z\left(\a\right)=0$ then $\a$ is perfect, complete and $\Der\g\cong\K\times\a$: the sequence terminates at $\K\times\a$, case 2. \\
If $Z\left(\a\right)\neq 0$ then we have 
\begin{equation}
	\Der\g\cong\ad\a\oplus\mathrm{I}
\end{equation}
(the center of $\ad\a$ is trivial when $\left[\a,\a\right]=\a$) with $\mathrm{I}=\K\epsilon\oplus\mathrm{I}_{12}$.\\ 
We have $\left[\epsilon,f\right]=-f$ for all $f\in\mathrm{I}_{12}$ such that the center of $\Der\g$ is zero and the sequence belongs to case 1.
\end{proof}
\begin{prop}\label{n17}Let $\a$ be a characteristic ideal of codimension 1 of $\g$. Then
\begin{enumerate}
\item $\dim\left(\Der\g\right)-\dim\left(\g\right)=\dim\left(\Der\g\mid_{\a}\right)-\dim\left(\ad\g\mid_{\a}\right)$.
	\item If $\dim\left(\Der\g\right)=\dim\left(\g\right)$, then 
		\item [ i)] $\Der\g\mid_\b=\ad\g\mid_\b$ for any ideal $\b$ contained in $\a$,
		\item [ ii)] $\g$ is algebraic,
		\item [ iii)] any ideal of codimension $1$ of $\g$ is characteristic.
\end{enumerate}
\end{prop}
\begin{proof} Let $\delta$ be a derivation of $\g$ vanishing on $\a$, then it vanishes on $\left[\g,\g\right]$ and $\left[\delta\cdot\g,\a\right]=0$. If $\a$ is not a direct factor of $\g$ $\delta\cdot\g\subset Z\left(\a\right)$. We check that all morphisms of $\g$ vanishing on $\a$ and with image contained in $Z\left(\a\right)$ are derivations of $\g$. Hence we have:
\begin{equation}\label{t1}
	\dim\left(\Der\g\mid_{\a}\right)=\dim\left(\Der\g\right)-\dim \left(Z\left(\a\right)\right)
\end{equation}
since 
\begin{equation}\label{t2}
	\dim\left(\g\right)=\dim\left(ad\g\mid_{\a}\right)+\dim \left(Z\left(\a\right)\right)
\end{equation}
from $\left(\ref{t1}\right)$ and $\left(\ref{t2}\right)$ the statement follows. If the ideal $\a$ is direct factor of $\g$ then it is perfect since $\a$ is characteristic and we directly verify the equality 1.\\
If $\dim\Der\g=\dim\g$ then 
\begin{equation}
	\Der\g\mid_\b=\ad\g\mid_\b
\end{equation} 
by 1. If the radical of $\g$ is not nilpotent then 
\begin{equation}
	\Der\g=\ad\g+\mathrm{J}
\end{equation}
where $\mathrm{J}$ is the ideal of derivations vanishing on $\b=\left[\g,\g\right]+\n$ and we have $\mathrm{J}^2=0$. Hence 
\begin{equation}
	\Der\g=\G\oplus N
\end{equation}
where $\G$ is a m.c.r Lie subalgebra of $e\left(\ad\g\right)$ and 
\begin{equation}
	N:=\mathrm{J}+ad\n
\end{equation}
is the largest nilpotent Lie subalgebra. We may check that 
\begin{equation}
	Z\left(\g\right)=Z\left(\n\right)^\G.
\end{equation} 
Let $\mathrm{I}$ be the ideal of morphisms of $\g$ which vanish on $\b$ and having their images contained in $Z\left(\g\right)$. Then 
\begin{equation}
	\mathrm{I}\cap\ad\g={0}
\end{equation}
indeed, if $x\in\g$ such that 
\begin{equation}
	\left[x,\g\right]\subset Z\left(\g\right)
\end{equation}
then $\G\cdot x=0$ and if $\left[x,\b\right]=0$ then $x\in Z\left(\g\right)$. Hence 
\begin{equation}
	\dim\left(\ad\g+\mathrm{I}\right)=\dim\left(\g\right)+\left(\mathrm{codim}\left(\b\right)-1\right)\cdot\dim \left(Z\left(\g\right)\right)
\end{equation}
If $\mathrm{codim}\left(\b\right)>1$ then the center of $\g$ is trivial and $\g$ is complete. Hence it is algebraic.
If $\mathrm{codim}\left(\b\right)=1$ then 
\begin{equation}
	\g=\s+\K\cdot x+n.
\end{equation}
By $\left(\ref{t1}\right)$, then 
\begin{equation}
	\Der\g=\ad\g\oplus\mathrm{I}
\end{equation}
and there exists an element $y\in\g$ such that $\ad y\mid{\a}$ is the semi-simple component of $\ad x\mid{\a}$. The decomposition $\s\oplus\K\cdot\ad y\oplus\n$ of $\g$ satisfies the following property
\begin{equation}
	\left[\s\oplus\K\cdot\ad y,\K\cdot\ad y\right]=0.
\end{equation} 
The torus $\K\cdot\ad y$ is maximal in the centralizer of $\ad\s$ in $\Der\g$ which is algebraic. Hence $\g$ is algebraic, cf. \cite{Ch}.
\end{proof}
\begin{cor}\label{n18}If $\dim\left(\Der\g\right)=\dim\left(\g\right)$ then the condition $Z\left(\g\right)\neq 0$ means that the codimension of $\left[\g,\g\right]$ in $\g$ is $1$ or $0$.
\end{cor}
\begin{thm}\label{n14}Let $\left(\g_n\right)_{n\in\N}$ be the derivation tower of a Lie algebra $\g$. Then it belongs to one of the following distinct cases:
\begin{enumerate}
	\item $Z\left(\g_n\right)=0$ for $n$ sufficiently large and the sequence terminates at a complete Lie algebra given by Theorem \ref{n10}.
	\item $Z\left(\g_n\right)\neq0$ for all $n$ and the sequence terminates at a Lie algebra $\gi$ equal to $\K\times\left[\gi,\gi\right]$ where $\left[\gi,\gi\right]$ is a perfect (i.e. equal to its derived ideal) and complete Lie algebra.
	\item The sequence of dimensions of $\g_n$ increases and diverges.
\end{enumerate}
\end{thm}
\begin{proof} 
First let us assume that if $(\g_n)$ does not satisfy 1., i.e. if $Z\left(\g_n\right)\neq 0$ for all $n$, then the sequence of dimensions $\dim\g_n$ increases in the large sense. If there would exist an integer $n$ such that $\dim\left(\Der\g\right)<\dim\g_n$
then $\g_n$ should be perfect, cf. \cite{Do}, and the center of $\Der\g_n$ should be zero since $\mathrm{I}=0$, a contradiction. Let us assume that $\left(\g_n\right)$ does not satisfy $3$., it just remains to study the sequences $\g_n$ such that $Z\left(\g_n\right)\neq 0$, with $\dim\g_n=\dim\g_p$ for all $n\geq p$ and we set $\g=\g_p$. From Corollary \ref{n18}, a Lie algebra satisfying 
\begin{equation}
	\dim\left(\Der\g\right)=\dim\g,\quad Z\left(\g\right)\neq 0,\quad\left[\g,\g\right]\neq\g
\end{equation}
admits an ideal $\a$ of codimension one equal to $\left[\g,\g\right]$. We may assume that $\left[\g,\g\right]\neq\g$ since a sequence associated with $\left[\g,\g\right]=\g$ belongs to case $1$. \\
If $Z\left(\g\right)$ is not included in $\a$, then we have 
\begin{equation}
	\g\cong\K\times\a
\end{equation}
and we may conclude in view of Lemma \ref{n16}. Suppose now that $Z(\g)\subset\a$. By Proposition \ref{n17}.2. we derive 
\begin{equation}
	\Der\g\mid_\a=\ad\g\mid_\a
\end{equation}
thus we may write: 
\begin{equation}
	\Der\g=\ad\g+\J
\end{equation}
where $\J$ is the ideal of derivations vanishing on $\a$. \\
If $\g=\K x\oplus\a$ for $x\notin\a$, $\J$ is the set of morphisms $\delta$ of $\g$ such that
\begin{equation}
	\delta\left(\a\right)=0\quad\mathrm{and}\quad\delta\left(x\right)\in Z_\g(\a)
\end{equation}
the centralizer of $\a$ in $\g$. We have 
$Z_\g\left(\a\right)\subset\a$ because 
\begin{equation}
	Z\left(\g\right)\subset\a,\quad\left[\Der\g,\Der\g\right]=\ad_\g\a
\end{equation}
since $\J^2=0$ and $\J\cdot\g\subset\a$. If $\Der\g$ has non-zero center, its derived ideal is also codimension 1 and the dimension of $\ad_\g\left(\a\right)$ is the same as $\a$, so $Z\left(\g\right)=0$, a contradiction.
\end{proof}

\begin{cor}\label{n19}
If the sequence $(\g_n)$ has an element $\g_p$ with non-nilpotent radical, then all $g_n$, $n\geq p$, have this property and the sequence is of type $1$ or $2$. The sequence of dimensions of type $3$ increases strictly from $\g_p$ on.
\end{cor}
\begin{proof} 
We will show that the sequence $\left(\g_n\right)$ associated with $\g$, with $\r\left(\g\right)\neq\n\left(\g\right)$, $Z\left(\g\right)\neq 0$ and $\dim\Der\g=\dim\g$ is of type $1$. The algebraic Lie algebra $\g$ has a decomposition $\s\oplus\u\oplus\n$ such that $\left[\s\oplus\u,\u\right]=0$ and the ideal $\left[\g,\g\right]=\s\oplus\n$ is of codimension $1$, by proof of theorem \ref{n14}. The ideal $C_2\g$ of $x\in\g$ such that $\left[x,\g\right]\subset Z\left(\g\right)$ is equal to $Z\left(\g\right)$ because we have $\left[x,\s\oplus\n\right]=0$ and that $\left[x,\u\right]\subset Z\left(\g\right)$ means $\left[x,\u\right]=0$, $\u$ consisting of semi-simple elements. So we have $\ad\g\cap \mathrm{I}=\ad\left(C_2\g\right)=0$ and $\Der\g$ is a direct product of the ideal $\ad\g$ with its center is trivial by the ideal $\mathrm{I}$ which satisfies $\mathrm{I}^2=0$. We conclude by using Lemma \ref{n15} showing that the center of $\Der\left(\Der\g\right)$ is trivial.
\end{proof}


\begin{thebibliography}{9}
\bibitem{Ch}C.~Chevalley, Theorie des alg\`ebres de Lie, Herman, 1968.
\bibitem{Do}J.~Dozias, Sur les d\'erivations des alg\`ebres de Lie. C. R. Acad. Sci. Paris Sér. I. Math. 259, p. 2748-2750 ( 1964).
\bibitem{Mo} G. D. Mostow, Fully reductive subalgebras of algebraic groups, Amer. J. Math. 68, p. 220-306 (1956).
\bibitem{Sc}E.~Schenkman, A theory of subinvariant Lie algebras. Amer. J. Math. 73, p. 453-474 (1951).
\end{thebibliography}
\end{document}